\documentclass[12pt]{elsarticle}
\usepackage{amsfonts}
\usepackage{amssymb}
\usepackage{mathrsfs}
\usepackage{amsmath}
\usepackage{amsthm}
\usepackage{hyperref}
\def\hD{\hat D}
\def\hH{\hat H}
\def\hX{\hat X}
\def\hc{\hat c}

\def\cL{\mathcal L}
\def\cG{\mathcal G}

\def\N{\mathbb{N}}

\def\tr{^{\mathsf{T}}}
\def\pr{^\prime}
\def\eop{\unskip\nobreak\hfil\penalty50\hskip2em\hbox{}\nobreak
\hfill\mbox{$\Box $}\par}
\newtheorem{thm}{Theorem}[section]
\newtheorem{prop}{Proposition}[section]
\newtheorem{lemma}{Lemma}[section]
\newtheorem{cor}{Corollary}[section]
\begin{document}
\begin{frontmatter}
	\title{Linearization and connection coefficients of polynomial sequences: A matrix approach }
	
\author{Luis Verde-Star}
 \address{
Department of Mathematics, Universidad Aut\'onoma Metropolitana, Iztapalapa,
Apartado 55-534, Mexico City 09340,
 Mexico }
\ead{verde@xanum.uam.mx}

 \begin{abstract}
	 For a sequence of polynomials  $\{p_k(t)\}_{k\ge 0}$ in one real or complex variable, where $p_k$ has degree $k$, for $k\ge 0$, we find explicit expressions and recurrence relations for infinite matrices whose entries are  the coefficients $d(n,m,k)$, called linearization coefficients, that satisfy
	 \begin{equation*}
p_n(t) p_m(t)=\sum_{k=0}^{n+m} d(n,m,k)\, p_k(t),\qquad n,m \in \N.
	 \end{equation*}
For any pair of polynomial  sequences  $\{u_k(t)\}_{k\ge 0}$ and  $\{p_k(t)\}_{k\ge 0}$ we find infinite matrices whose entries are the  coefficients $e(n,m,k)$ that satisfy
	 \begin{equation*}
p_n(t) p_m(t)=\sum_{k=0}^{n+m} e(n,m,k)\, u_k(t),\qquad n,m \in \N.
	 \end{equation*}
Such results are obtained using the matrix approach of our previous papers \cite{Opm} and \cite{PSGH}.
We also obtain recurrence relations for the linearization coefficients,  apply the general results to general orthogonal polynomial sequences and to  particular families of orthogonal polynomials  such as the Chebyshev, Hermite, and Charlier families.



\end{abstract}
\end{frontmatter}

\section{Introduction}
A sequence of polynomials $\{p_n(t)\}_{n\ge 0}$ in one real or complex variable such that $p_n$ has degree $n$, for $n \ge 0$, is called a polynomial sequence and it is a basis for the algebra of all polynomials in $t$. Therefore every polynomial $u_r(t)$ of degree $r$ has a unique representation of the form
\begin{equation*}
u_r(t)= \sum_{k=0}^r d_k p_k(t),
\end{equation*}
and hence  finding the coefficients $d_k$ for a given $u_m$ is a computational problem that appears often in many areas of Mathematics and its applications.
An important particular case is obtained when $u_r(t)=p_n(t) p_m(t)$. In this case we have
\begin{equation*}
p_n(t) p_m(t)=\sum_{k=0}^{n+m} d(n,m,k) p_k(t), \qquad n,m \in \N,
\end{equation*}
and the coefficients $d(n,m,k)$ are called the {\em linearization coefficients} of the sequence $\{p_k(t)\}_{k\ge 0}$. A related problem is obtained when we want to find the coefficients $e(n,m,k)$ in 
\begin{equation*}
p_n(t) p_m(t) = \sum_{k=0}^{m+n} e(n,m,k) u_k(t), \qquad n,m \in \N,
\end{equation*}
where $\{u_k(t)\}_{k\ge 0}$ is another polynomial sequence. 

In this paper, for any given pair of polynomial sequences $\{p_k(t)\}_{k\ge}$ and $\{u_k(t)\}_{k\ge 0}$ 
we find explicit expressions and recurrence relations for  certain matrices whose entries are  the linearization coefficients $d(n,m,k)$ and $e(n,m,k)$. This is accomplished using some results from \cite{PSGH} about Hessenberg matrices, and the matrix approach that we have used in \cite{Opm} and \cite{Uni}. If the polynomial sequences have certain properties then it is possible to simplify the explicit expressions and the recurrence relations for the matrices of linearization coefficients. This happens for orthogonal polynomial sequences, sequences of binomial type, Sheffer sequences, sequences of interpolatory type, and polynomial sequences associated with infinite matrices that are representations of linear operators on the algebra of formal Laurent series in one variable, such as multiplication,   composition, and differential  operators. There is an interesting matrix approach to study several classes of polynomial sequences in the recent book \cite{Costa}.

In this paper, in addition to the case of general polynomial sequences  we consider only the class of general orthogonal polynomial sequences and give examples of our results for the  Chebyshev, Hermite, and Charlier  families of  polynomial sequences.

In the literature there are numerous papers about linearization coefficients of particular polynomial sequences. Most of them deal with orthogonal polynomials and use analytical methods. See \cite[Ch. 9]{Ismail} and the  references therein.

Finding connection coefficients is essentially a problem of change of bases in the space of polynomials. For polynomial sequences in certain classes, such as the Sheffer sequences, the connection coefficients are easily  obtained using the methods of Umbral calculus. See \cite{Roman} and \cite{Rota}. 

In Section 2 we present some basic material about infinite matrices and a result from \cite{PSGH} that describes the connection between Hessenberg matrices and invertible lower triangular matrices.
In Section 3 we obtain the general results for linearization and connection coefficients of general polynomial sequences. In Section 4 we apply the general results to general orthogonal polynomial sequences and in Section 5 we consider the particular families of Chebyshev, Hermite, and Charlier orthogonal polynomials.

The application of our general results to other classes of polynomial sequences will be presented in a forthcoming paper. 

  \section{Infinite  generalized lower Hessenberg matrices}
In this section we define the algebra $\cL$  of infinite generalized lower Hessenberg matrices and present some of their basic properties  that will be  used in the rest of the paper. 
 For additional information about the algebra $\cL$ see \cite{Opm} and \cite{PSGH}.

An infinite matrix  $A=[a_{j,k}]$, where the indices run over the non-negative integers and the entries are complex  numbers is  a {\sl lower generalized Hessenberg matrix}
 if  there exists an integer $m$ such that  $a_{j,k}=0$ whenever
 $j-k < m$. We denote by $\cL$ the set of all such matrices.
 
 We say that the entry $a_{j,k}$ of $A$ lies in the $n$-th diagonal of $A$ if $j-k=n$.  If $m > n$ then the $m$-th diagonal lies below (to the left of) the $n$-th diagonal.
 A nonzero element of $\cL$ is a {\sl diagonal matrix} if all of its nonzero elements lie in a single diagonal. 

 If  $A$ is  a nonzero element of $\cL$ and  $m$ is  the minimum integer such that $A$ has at least one nonzero entry in the $m$-th diagonal, then we say that $A$ has {\sl index} $m$ and write ind$(A)=m$. The index of the zero matrix is infinity, by definition.
 It is clear   that $\cL$ is a complex vector space with the natural addition of matrices and multiplication by scalars. It is also closed under matrix multiplication.
If $A$ and $B$ are in $\cL$, with ind$(A)=m$ and ind$(B)=n$, then the product $C=AB$ is a well defined element of $\cL$  and
\begin{equation}\label{eq:mult}
	c_{i,k}= \sum_{j=k+n}^{i-m} a_{i,j} b_{j,k}, \qquad i-k \ge m+n.
\end{equation} 
Note that ind$(AB)\ge m +n$ and that the multiplication in $\cL$ involves only finite sums.

 A  sufficient, but not necessary,  condition for $A$ to have a two-sided inverse is that ind$(A)=0$ and  $a_{k,k}\ne 0$ for $k \ge 0$. We denote by $\cG$ the set of all matrices that satisfy such condition. It is clear that $\cG$ is a group under matrix multiplication. The unit is the identity matrix $I$ whose entries on the 0-th diagonal are equal to 1 and all other entries are zero.

We  define next  some particular elements of $\cL$ that will be used often in the rest of the paper. Let  $X$  denote the diagonal matrix of index -1 with $X_{j,j+1}=1$ for $j \ge 0$, and denote by $\hat X$ the transpose of $X$. Note that $\hat X$ is diagonal of index  1,   $X \hat X= I$ and $\hat X X= J_0$, where $J_0$ is the diagonal matrix  of index zero that has its entry in the position $(0,0)$ equal to zero and its entries in positions $(j,j)$ equal to 1 for $j \ge 1$.
Therefore, $\hat X$ is a right-inverse for $X$, but it is  not a left-inverse.

We say that a matrix  $A \in \cL$ of index $m$ is {\sl monic} if all the entries in the  diagonal of index $m$ are equal to 1. Note that  a monic matrix of index -1 is a unit lower  Hessenberg matrix.

If $m$ is a positive integer then $X^m$ is the diagonal matrix of index $-m$ with all its entries in the $-m$ diagonal equal to 1. Analogously, $\hX^m$ is diagonal of index $m$ and all its entries in the $m$-th diagonal are equal to 1.

The following theorem, which  was proved in \cite{PSGH}, will be used in Section 3 to obtain our results about linearization and connection coefficients.

\begin{thm}\label{simil}
Let $H$ be a monic matrix of index -1. Then there exists a unique monic  $A$  in $ \cG$ such that $ A H = X A $.
\end{thm}
{\it Proof:}  For $j \ge 0 $ let $r_j$ denote the $j$-th row of the identity matrix $I$. Then, for any matrix $Z$ the product $r_j Z$ is the $j$-th row of $Z$. In particular, $r_j X = r_{j+1}$. Therefore, if $A$ satisfies $ A H = X A $ then we have  $ r_j A H = r_j X A = r_{j+1} A $ for $j \ge 0$.
This means that if we know the $j$-th row of $A$ then the $(j+1)$-th row is obtained by multiplying the $j$-th row by $H$ on the right. If we take the $0$-th row of $A$ equal to $r_0$ then we can construct $A$ row by row by repeated multiplication by $H$. It is clear that the resulting matrix $A$ is monic of index zero and therefore it is in the group $\cG$. The uniqueness of $A$ is also clear, since every monic element of $\cG$ has its $0$-th row equal to $r_0$.   \eop

The previous theorem says that every monic matrix of index -1 is similar to $X$.
 The proof of the theorem can be easily modified to show that every monic matrix of index $-m$  is similar to $X^m$, for $ m >1$, but in that case $A$ is not unique.

It is easy to see  that for every monic matrix  $A$  in $ \cG$ the matrix  $A^{-1} X A$ is  monic  of index -1.

\section{Polynomial sequences}
A  matrix $A=[a_{k,j}]$ of  index $m $  determines a  sequence of polynomials  $u_k(t)$ defined by 
\begin{equation}\label{eq:rowpoly}
u_k(t) = \sum_{j=0}^{k-m}  a_{k,j} t^j, \qquad k \ge 0.
\end{equation}
That is, the entries in the $k$-th row of $A$ are the coefficients of $u_k$. 
If $m> 0$ then $u_j =0$  for $ 0 \le j \le m-1$, and if $A$ is monic then each nonzero  $u_k$ is monic and has degree $k-m$.

If $v_k(t)$ are the polynomials associated with $X A$ then  $ v_k(t)= u_{k+1}(t)$ for $ k \ge 0$. If $ w_k(t)$ are the polynomials associated with $A X$ then $w_k(t) = t u_k(t)$ for $k \ge 0$.
 
Define $C_t=[1,t,t^2,\ldots ]\tr$, where $t$ is a real or complex variable.
Then for each $A$ in $\cG$ we have
\begin{equation} \label{ACt}
 A C_t=[u_0(t), u_1(t), u_2(t),\ldots]\tr,
\end{equation}
 $\{u_n(t)\}_{n\ge 0}$ is a polynomial sequence and it is monic if $A$ is monic.

\begin{prop}\label{multOp}
Let $A$ be in $\cG$, let $w(t)$ be a polynomial, and $C_t$ and  $\{u_n(t)\}_{n\ge 0}$  as defined above. Then we have
	\begin{equation*}
	A \, w(X) \, C_t=[u_0(t) w(t), u_1(t) w(t),u_2(t) w(t),\ldots]\tr.
	\end{equation*}
\end{prop}

{\em Proof.} It follows by linearity from
	\begin{equation*}
	A \, X  \, C_t=[t u_0(t),t  u_1(t),t u_2(t),\ldots]\tr.
\end{equation*}
\eop

From Theorem \ref{simil} we see that each monic  matrix $H=[h_{k,j}]$ of index -1 is associated with two invertible monic  matrices $A$ and $P=A^{-1}$ of index zero that satisfy $ A H P = X $. Let us denote by $u_k(t)$ and $p_k(t)$ the polynomial sequences determined by $A$ and $P$, respectively. 

The matrix equation $A H = X A$  says that the $j$-th row of $A$  multiplied on the right by $H$ gives us  the $(j+1)$-th row of $A$. This fact was used in the proof of Theorem \ref{simil} and can be used to construct the sequence of polynomials $u_k(t)$.  

From the  equation $ H P = P X$ we get immediately
\begin{equation}\label{eq:recrel}
  p_{k+1}(t) +   \sum_{j=0}^k h_{k,j} p_j(t) = t p_k(t), \qquad k\ge 0,
\end{equation} 
which is a recurrence relation that can be used to compute $p_{k+1}(t)$ as a linear combination of $ p_0(t),p_1(t),\ldots ,p_k(t)$ and $t p_k(t)$.
When  $H$ is tridiagonal \eqref{eq:recrel} is  the well-known three-term recurrence relation for orthogonal polynomial sequences. In that case the $p_k$ are orthogonal with respect to some linear functional on the space of polynomials. If  $H$ is tridiagonal then the equation  $A H = X A$ yields a three-term recurrence relation  for the columns of $A$.  

The matrix $P=A^{-1}$ can  be constructed in a way analogous to the construction of $A$ in the previous section, but instead of using $H$ as a multiplier in each step, we use a one-sided inverse  $\hH$  of $H$, and $P$ is computed column by column. This is done as follows.
Let us define  $Y= H  \hX$ and   $\hH= \hX Y^{-1}$. Note that  $Y$ is a monic element of $ \cG$ and  $\hH$ is monic of index 1. It is easy to verify that $ H  \hH =I$, that is $\hH$ is a right inverse for $H$, and  $\hH  H$ differs from the identity matrix only in the 0-th column. Since $X \hX=I$, from the equation $H P = P X$ we get $ \hH H P \hX = \hH P$. Since $P$ is lower triangular, $\hH H P$ differs from $P$ only in the (0,0) entry. Therefore $P \hX = \hH P$ holds.  
This means that $\hH$ times the $k$-th column of $P$ equals the $(k+1)$-th column of $P$. 
Let us note  that the 0-th column of $P$ is equal to the 0-th column of $- \hH H$, with its 0-th entry set equal to 1. Therefore $P$ can be computed column by column by repeated multiplication by $\hH$. 

The matrix $P$  can also be  constructed  using the recurrence relation \eqref{eq:recrel}, as it is usually done in the case of orthogonal polynomial sequences, but note that in the general case  $p_k$ depends on $p_0,p_1,\ldots, p_{k-1}$.

\begin{thm}\label{genlinear}
	Let $P$ be a monic element of $\cG$ and let   $H= P X P^{-1}$. Let   $\{p_k(t)\}$ be the polynomial sequence associated with $P$. Then, for every   polynomial $w$ of degree $m$ we have
	\begin{equation}\label{eq:wofH}
	P\, w(X) \, C_t =w(H)\, P \, C_t,
	\end{equation}
	and 
	\begin{equation}\label{eq:wlinear}
		p_n(t) w(t)=\sum_{k=0}^{n+m} w(H)_{n,k}\, p_k(t), \qquad n \in \N.
	\end{equation}
\end{thm}

{\em Proof.} Since $X=P^{-1} H P$ we have $w(X)=P^{-1} w(H) P$ and thus $P\, w(X)=w(H)\, P$. Multiplying both sides of  this equation by the column vector $C_t$ we obtain \eqref{eq:wofH}. Comparing the $n$-th rows in both sides of \eqref{eq:wofH} and using Proposition \ref{multOp} we obtain \eqref{eq:wlinear}.
\eop

Taking $w=p_m$ we obtain immediately the following corollary.

\begin{cor}\label{plinear}
With the hypothesis of the previous theorem we have
	\begin{equation}\label{eq:pnmkLin}
p_n(t) p_m(t)= \sum_{k=0}^{n+m} d(n,m,k) p_k(t), \qquad n,m \in \N,
	\end{equation}
where  
	\begin{equation}\label{eq:dnmk}
	d(n,m,k)=p_m(H)_{n,k}, \qquad n,m,k \in \N.
	\end{equation}
\end{cor}

The numbers $d(n,m,k)$ are the linearization coefficients of the polynomial sequence $\{p_k(t)\}_{k\ge 0}$.

From $ n,m,k$ in $\N$ we obtain  from  equation  \eqref{eq:pnmkLin} the following properties
\begin{enumerate}
	\item{ } $d(n,m,k)=d(m,n,k)$,

	\item{ } $d(n,m,k)=0$ if $n+m<k$,

	\item{ }  $d(n,m,k)=1$ if $n+m=k$,

	\item{ }  $d(0,m,k)=\delta_{m,k}$.
\end{enumerate}
Let us note that \eqref{eq:dnmk} and  $d(n,m,k)=d(m,n,k)$ give us 
\begin{equation}\label{eq:rows}
\hbox{Row}_m(p_n(H))=\hbox{Row}_n(p_m(H)), \qquad n,m \ge 0,
\end{equation}
and thus the rows of indices $ 0,1,2,\ldots, n-1$ of $p_n(H)$  appear as the $n$-th row of 
$p_0(H), p_1(H), \ldots p_{n-1}(H).$

\begin{thm}\label{recLinCoef}
The linearization coefficients $d(n,m,k)$ satisfy the recurrence relation
\begin{equation}\label{eq:recdnmk}
	d(n+1,m,k)=d(n,m+1,k) +(h_{m,m} - h_{n,n}) d(n,m,k) +  \sum_{j=0}^{m-1}  h_{m,j} d(n,j,k) - \sum_{j=0}^{n-1} h_{n,j} d(m,j,k).
\end{equation}
\end{thm}

{\em Proof.} The recurrence relation \eqref{eq:recrel} gives us the matrix equations
\begin{equation}\label{eq:recpnH}
p_{n+1}(H)= H p_n(H) - \sum_{j=0}^n h_{n,j} p_j(H),\qquad n \ge 0.
\end{equation}
 For the entries with indices $(m,k)$ of the matrices in the previous equations we obtain 
 \begin{equation*}
 p_{n+1}(H)_{m,k} =( H p_n(H))_{m,k} - \sum_{j=0}^n h_{n,j} p_j(H)_{m,k},
 \end{equation*}
 and using the definition of the linearization coefficients $d(n,m,k)$  we get
 \begin{equation*}
	 d(n+1,m,k)= \sum_{j=0}^{m+1} h_{m,j} d(n,j,k)-\sum_{j=0}^n h_{n,j} d(j,m,k).
 \end{equation*}
Since $h_{m,m+1}=1$, for $m \ge 0$, and  $d(n,m,k)=d(m,n,k)$, this recurrence relation  can be written as 
\begin{equation*}
	d(n+1,m,k)=d(n,m+1,k) +(h_{m,m} - h_{n,n}) d(n,m,k) +  \sum_{j=0}^{m-1}  h_{m,j} d(n,j,k) - \sum_{j=0}^{n-1}  h_{n,j} d(m,j,k),
\end{equation*}
 which is \eqref{eq:recdnmk}.   \eop

 From the recurrence relation \eqref{eq:recpnH} we obtain immediately
 \begin{equation}\label{eq:rowspnH}
	 \hbox{Row}_m(p_{n+1}(H))= \sum_{j=0}^{m+1} h_{m,j} \hbox{Row}_j(p_n(H)) - \sum_{j=0}^n h_{n,j} \hbox{Row}_m(p_j(H)), \qquad n,m \ge 0.
 \end{equation}

 For fixed $k$  the recurrence relation \eqref{eq:recdnmk} can be used to compute the coefficients $d(n+1,m,k)$, for $m\ge 0$ if we know the coefficients $d(n,m,k)$, for $m \ge 0$. Since $d(0,m,k)=\delta_{m,k}$, for fixed $k$ we can compute the matrix $d(n,m,k)$ row by row. Recall that $d(n,m,k)=d(m,n,k)$.

 \begin{thm}\label{linpu}
	 Let $P$ and $U$ be monic elements of $\cG$ with associated polynomial sequences $\{p_k(t)\}_{k\ge 0}$ and $\{u_k(t)\}_{k\ge 0}$ respectively. Define $H=P  X  P^{-1} $ and $ K=U  X  U^{-1}$.  Then we have
\begin{equation}\label{eq:linearpu}
p_n(t) p_m(t)= \sum_{k=0}^{n+m} e(n,m,k) u_k(t), \qquad n,m \in \N,
\end{equation}
where
	 \begin{equation}\label{eq:enmk}
	e(n,m,k)= \sum_{j=0}^{n+m} p_n(H)_{m,j} p_j(K)_{0,k}, \qquad n,m,k \in \N.
\end{equation}
 \end{thm}
 {\em Proof.} For $j \ge 0$ let $r_j$ denote the $j$-th row of the identity matrix.  Since  $ X=U^{-1} K U$ we have $ p_j(X)= U^{-1} p_j(K) U$. Therefore
\begin{equation*}
	r_j P = r_0 p_j(X)=r_0 U^{-1} p_j(K) U= r_0 p_j(K) U,
\end{equation*} 
and hence
\begin{equation*}
r_m(p_n(H) P) = \left(\sum_{j=0}^{n+m} p_n(H)_{m,j} \, r_0 p_j(K)\right) U.
\end{equation*}
Multiplying both sides of this equation by the column vector  $C_t$ we obtain \eqref{eq:linearpu}. \eop

Taking $n=0$ in the previous theorem we obtain the following corollary.

\begin{cor}
With the hypothesis of the previous theorem, for $m \ge 0$ we have
	\begin{equation}\label{eq:connectpu}
p_m(t)= \sum_{k=0}^m e(0,m,k) u_k(t)=\sum_{k=0}^m p_m(K)_{0,k}\, u_k(t).
	\end{equation}
\end{cor}

The numbers $p_m(K)_{0,k}$ are the connection coefficients of the sequences $\{p_n(t)\}_{n\ge 0}$ and $\{u_n(t)\}_{n\ge 0}$.

It is clear that the inverse relation of \eqref{eq:connectpu} is 
	\begin{equation}\label{eq:connectup}
u_m(t)= \sum_{k=0}^m u_m(H)_{0,k}\,  p_k(t),
	\end{equation}
and therefore
\begin{equation}\label{eq:inverseConn}
	\sum_{k=0}^{m+n} p_m(K)_{0,k} \, u_k(H)_{0,n}=\delta_{m,n}.
\end{equation}

\section{Orthogonal polynomial sequences.}
In this section we consider monic matrices $P$  in $\cG$  for which the Hessenberg matrix $H=P X P^{-1}$ is tridiagonal of index -1  and has nonzero elements in its diagonal of index 1. We will show that the associated polynomial sequences of such matrices are orthogonal with respect to a linear functional defined on the space of polynomials.  We will apply the general results of the previous sections to the orthogonal polynomial sequences.

Let $H$ be a monic tridiagonal matrix of index -1 with entries $H_{k,k}=\beta_k$,  and $H_{k+1,k}=\alpha_{k+1}$, for $k \ge 0$, and  such that $\alpha_k \ne 0$ for $k\ge 1.$ Let $P=[c_{n,k}]$ be the unique monic element of $\cG$ that satisfies $H P= P X$, and let $\{p_k(t)\}_{k\ge 0}$ be the  polynomial sequence associated with $P$. Let  $P^{-1}=[\hc_{n,k}]$. The entries of $P$ and $P^{-1}$ can be expressed in terms of the entries of $H$. See \cite{Uni}.

 In this case the recurrence relation \eqref{eq:recrel} becomes the three-term recurrence 
 \begin{equation}\label{eq:3term}
p_{n+1}(t)= (t- \beta_n) p_n(t) - \alpha_k p_{n-1}(t), \qquad n\ge 1.
 \end{equation}

The recurrence relation \eqref{eq:recdnmk}  for the linearization coefficients $d(n,m,k)=p_n(H)_{m,k}$ becomes 
\begin{equation}\label{eq:rrdOP}
	d(n+1,m,k)=d(n,m+1,k) +(\beta_m -\beta_n) d(n,m,k) +\alpha_m d(n,m-1,k) -\alpha_n d(n-1,m,k).
\end{equation}
Let us recall that $d(n,m,k)=d(m,n,k)$ for all $n,m,k$ in $\N$.

\begin{lemma}\label{dnm0}
	The linearization coefficients $d(n,m,k)$ of the  sequence $\{p_n(t)\}_{n\ge 0}$ satisfy
\begin{equation}\label{eq:dnmoDiag}
d(n,m,0)= \delta_{n,m} \alpha_1 \alpha_2 \cdots \alpha_n, \qquad n,m \in \N.
\end{equation}
\end{lemma}
{\em Proof.} From the linearization equation
\begin{equation*}
p_n(t) p_m(t)=\sum_{k=0}^{n+m} d(n,m,k) p_k(t)
\end{equation*}
we obtain  $p_0(t) p_m(t)= d(0,m,m) p_m(t)$ and thus $d(0,0,0)=1$ and $d(0,m,0)=0$ for $m\ge 1$.

Since $p_1(t)= t-\beta_0$ we have
\begin{equation*}
	p_1(t) p_m(t) = (t-\beta_m+\beta_m -\beta_0) p_m(t) -\alpha_m p_{m-1}(t)+\alpha_m  p_{m-1}(t).
\end{equation*}
By the three-term recurrence relation we get
\begin{equation*}
	p_1(t) p_m(t) = p_{m+1}(t)  + \alpha_m p_{m-1}(t) +(\beta_m -\beta_0) p_m(t).
\end{equation*}
Since  $\{p_k, k \in \N\}$ is a basis for the vector space of polynomials the previous equation gives us $d(1,m,0)=\delta_{1,m} \alpha_1$, for $m \ge 0$. Therefore the first two rows of the symmetric matrix with entries $d(n,m,0)$ are
\begin{equation*}
	\left[ \begin{matrix} 1 & 0 & 0 & 0 & 0 &  \cdots \cr
		0 & \alpha_1 & 0 & 0 & 0 & \cdots \cr
	\end{matrix} \right].
\end{equation*}
Using the recurrence relation \eqref{eq:rrdOP} we see that the only nonzero entries appear in the main diagonal and the recurrence becomes $d(n,n,0) = \alpha_n d(n-1,n-1,0)$. Therefore $d(n,m,0)=\delta_{n,m} \alpha_1 \alpha_2 \cdots \alpha_n.$  \eop

We define the linear functional $\tau$ on the space of polynomials by
\begin{equation}\label{eq:tau}
\tau(t^k)= \hc_{k,0}, \qquad k\in \N.
\end{equation}
	Lets recall that $\hc_{n,k}$ are the entries of $P^{-1}$.

\begin{thm}\label{orthogonal}
The polynomial sequence $\{p_n(t)\}_{n\ge 0}$ is orthogonal with respect to the linear functional $\tau$, that is, for $n$ and $m$ in $\N$ we have
	\begin{equation}\label{eq:orthopntau}
\tau(p_n(t) p_m(t))= \delta_{n,m} \alpha_1 \alpha_2 \cdots \alpha_n.
	\end{equation}
\end{thm}
{\em Proof.} Since $H P=P X$ we have $p_n(H)P= P p_n(X)$. By Proposition \ref{multOp} the entries in the  $m$-th row of $P p_n(X)$  are the coefficients of $p_n(t) p_m(t)$ with respect to the basis of monomials, that is,
\begin{equation*}
	p_n(t) p_m(t)=\sum_{j=0}^{m+n} (P p_n(X))_{m,j} t^j.
\end{equation*}
Therefore 
\begin{equation*}
	(P p_n(X) P^{-1})_{m,0}=\sum_{j=0}^{m+n} (P p_n(X))_{m,j} P^{-1}_{j,0}= \sum_{j=0}^{m+n} (P p_n(X))_{m,j} \tau(t^j),
	\end{equation*}
and by the definition of the coefficients $d(n,m,k)$ and Lemma \ref{dnm0} this equation is equivalent to
\begin{equation*}
	\tau(p_n(t) p_m(t))= p_n(H)_{m,0} = d(n,m,0)=\delta_{n,m}\alpha_1 \alpha_2 \cdots \alpha_n.
\end{equation*}
\eop

Using the recurrence relation for the linearization coefficients $d(n,m,k)$ we can show that $d(n,m,k)=0$ if $ k <|n-m|$ and therefore 
\begin{equation}\label{eq:linOPsimple}
p_n(t) p_m(t)=\sum_{k=|n-m|}^{n+m} d(n,m,k) p_k(t).
\end{equation}

{\em Remark.}  If instead of a tridiagonal matrix $H$ we consider a banded monic Hessenberg matrix of order -1,  with four diagonals and with nonzero entries in the diagonal of index 2, then we can show that the corresponding polynomial sequence  $\{p_n(t)\}_{n\ge 0}$ is partially  orthogonal, that is, for $n\ge 0$  the polynomial $p_n$ is orthogonal to $p_m$ if $m \ge 2n +1$. When $H$ is pentadiagonal then $p_n$ is orthogonal to $p_m$ if $m \ge 3 n +1.$  

\section{Some families of orthogonal polynomial sequences}

In this section we compute the linearization coefficients of some simple families of orthogonal polynomial sequences. For each family  we find  explicit expressions for the matrices $p_n(H)$ in terms of the coefficients of the three-term recurrence relation.

 Let  $D$ be  the diagonal matrix of index 1 that has $D_{k+1,k}=k+1$  for $k \ge 0$, and all the other entries equal to zero. 
If the $k$-th row of a matrix $P$ corresponds to a polynomial $p_k(t)$ then the $k$-th row of $PD$ corresponds to $p_k\pr (t)$. The matrix $D$ satisfies $X D - D X=I$, where $I$ is the identity matrix. For $k \ge 0$ the index of $D^k$ is $k$. 

We consider first the Chebyshev family of orthogonal polynomial sequences. 
 Let $a$ and $b$ be complex numbers with $a \ne 0$ and let $H=a \hX + b I + X$.
Note that $H$ is a monic Toeplitz tridiagonal matrix of index -1.
Let $P$ be the unique monic element of $\cG$ that satisfies $ H P = P X$ and let $\{p_n(t)\}_{n\ge 0}$ be the polynomial sequence associated with $P$.
If $a=1/4$ and $b=0$ then $\{p_n(t)\}_{n\ge 0}$ is the sequence of monic Chebyshev polynomials of the first kind. In \cite{PSGH} we studied the polynomial sequences associated with general Hessenberg-Toeplitz matrices.

The generating function of the $k$-column of $P$ is $z^k/(1+b z +a z^2)^{k+1}$ and therefore $P$ is obtained by deleting the $0$-th row and the $0$-th column  from the matrix that represents the composition operator on the vector space of formal power series that sends  $z^k$  to $ z^k/(1+b z +a z^2)^{k}$.
 We also have the series representation  
\begin{equation}\label{eq:PCheby}
	P=\sum_{k=0}^\infty  ( X-H)^k\, \dfrac{D^k}{k!}.
\end{equation}
It is not an exponential series, since $X-H$ and $D$ do not commute. 

The inverse of $P$ is also a truncated composition matrix that taking   $a=1$ and $b=2$ becomes the Catalan triangle.

Using the recurrence relation \eqref{eq:recpnH} it is easy to show that in this case 
\begin{equation}\label{eq:pnHCheby}
p_n(H)=\sum_{k=0}^n a^k \hX^k \, X^{n-k}, \qquad n\ge 0,
\end{equation}
where $\hX$ is the transpose of $X$. Note that $p_n(H)$ is independent of $b$.
For example
\begin{equation*}
	p_3(H)= \left[ \begin{matrix} 0 & 0 &0 & 1 & 0 & 0 & \ldots \cr
		0 & 0 & a & 0 & 1 & 0 & \ldots \cr
			   0 & a^2 & 0 & a & 0 & 1 & \ldots  \cr
			   a^3 & 0 & a^2 & 0 & a & 0 & \ldots  \cr
			   0 & a^3 & 0 & a^2 & 0 & a & \ldots  \cr
			  0& 0 & a^3 & 0 & a^2 & 0  & \ldots  \cr
	\vdots & \vdots & \vdots & \vdots & \vdots & \vdots &\ddots  \cr \end{matrix} \right].
\end{equation*}

The matrix $A$ with $A_{n,m}=d(n,m,0)$ is the diagonal matrix of index zero  $A=\hbox{Diag}(1,a,a^2,\ldots)$, and for $k\ge 0$ the matrix $A^{(k)}$ with $(A^{(k)})_{n,m}=d(n,m,k)$ has the representation 
\begin{equation}\label{eq:AknmCheby}
	A^{(k)}=\sum_{j=0}^k  \hX^{j}\,  A \, X^{k-j}, \qquad k \ge 0.
\end{equation}
For example
\begin{equation*}
	A^{(3)}= \left[ \begin{matrix} 0 & 0 &0 & 1 & 0 & 0 & 0 &\ldots \cr
		0 & 0 & 1 & 0 & a & 0 & 0 &  \ldots \cr
			   0 & 1 & 0 & a & 0 & a^2 & 0 & \ldots  \cr
			   1 & 0 & a & 0 & a^2 & 0 &  a^3 &\ldots  \cr
			   0 & a & 0 & a^2 & 0 & a^3 & 0 & \ldots  \cr
			  0 & 0 & a^2 & 0 & a^3 & 0  & a^4 & \ldots  \cr
		0&  0& 0 & a^3 & 0 & a^4 & 0   & \ldots  \cr
	\vdots & \vdots & \vdots & \vdots & \vdots & \vdots & \vdots  &\ddots  \cr \end{matrix} \right].
\end{equation*}
Let us note that $A^{(k)}$ is a symmetric matrix of index $-k$.

We consider next a family of generalized Hermite polynomials with two parameters.
Let $a$ and $b$ be complex numbers with $a \ne 0$ and define 
\begin{equation}\label{eq:HHermite}
H=X + b I + a D.
\end{equation}
 The unique monic matrix  $P$  in $\cG$ that satisfies $H P= P X$ has the exponential series representation
 \begin{equation}\label{eq:HermExp}
	 P=\sum_{k=0}^\infty \dfrac{1}{k!} \left(-b D - a \dfrac{D^2}{2}\right)^k,
\end{equation}
and therefore 
 \begin{equation}\label{eq:InvHermExp}
	 P^{-1}=\sum_{k=0}^\infty \dfrac{1}{k!} \left(b D + a \dfrac{D^2}{2}\right)^k.
\end{equation}

Let $\{p_n(t)\}_{n\ge 0}$ be the polynomial sequence associated with $P$. From \eqref{eq:HermExp} and \eqref{eq:InvHermExp} we can see that the  moments, that is, the entries in the $0$-th column of $P^{-1}$, are obtained from the $0$-th column of $P$ replacing $a$ with $-a$ and $b$ with $-b$.

Using induction and the three-term recurrence relation we can show that
\begin{equation}\label{eq:pnHHermite}
p_n(H)=\sum_{k=0}^n \binom{n}{k} a^k D^k \, X^{n-k}, \qquad 
 n\ge 0.
\end{equation}
Note that $p_n(H)$ is independent of the parameter $b$.

For the Hermite family  the matrix $A$ with $A_{n,m}=d(n,m,0)$ is the diagonal matrix of index zero  $A=\hbox{Diag}(1,a,2! a^2,3! a^3, \ldots)$, and for $k\ge 0$ the matrix $A^{(k)}$ with $(A^{(k)})_{n,m}=d(n,m,k)$ has the representation 
\begin{equation}\label{eq:AknmHerm}
	A^{(k)}= \dfrac{1}{k!} \sum_{j=0}^k \binom{k}{j}  D^j   A \hD^{k-j}, \qquad k \ge 0.
\end{equation}
where $\hD$ is the transpose of $D$.

Our last example is the family of Charlier polynomials.
Let $a$ be a nonzero complex number and define
\begin{equation}\label{eq:HCharlier}
H= X + X D +(a-1) I + a D.
\end{equation}
Let $P$ be the unique monic element of $\cG$ that satisfies $ H P = P X$ and let $\{p_n(t)\}_{n\ge 0}$ be the polynomial sequence associated with $P$.

In this case the matrices $p_n(H)$ can be expressed as
\begin{equation}\label{eq:pnHCharlier}
	p_n(H)= \sum_{k=0}^n \binom{n}{k} a^k D^k (I +D)^{n-k} X^{n-k}, \qquad n \ge 0.
\end{equation}

For the Charlier  family  the matrix $A$ with $A_{n,m}=d(n,m,0)$ is the diagonal matrix of index zero  $A=\hbox{Diag}(1,a,2! a^2,3! a^3, \ldots)$, and for $k\ge 0$ the matrix $A^{(k)}$ with $(A^{(k)})_{n,m}=d(n,m,k)$ has the representation 
\begin{equation}\label{eq:AknmCh}
A^{(k)}=\dfrac{1}{k!}  \sum_{j=0}^k \binom{k}{j}  D^j   (I +D)^{k-j}    A  \hD^{k-j}, \qquad k\ge 0,
\end{equation}
where $\hD$ is the transpose of $D$.

\end{document}